\documentclass[twoside]{amsart}
\usepackage{amsmath}
\usepackage{amsfonts}
\usepackage{amssymb,enumerate}
\usepackage{amsthm}
\usepackage{hyperref}

\theoremstyle{plain}
\newtheorem{theorem}{Theorem}[section]
\newtheorem{lemma}[theorem]{Lemma}
\newtheorem{definition}[theorem]{Definition}
\newtheorem{corollary}[theorem]{Corollary}
\newtheorem{proposition}[theorem]{Proposition}
\newtheorem{remark}[theorem]{Remark}

\newtheorem{question}[theorem]{Question}
\begin{document}

\title{ Model Theory of the Inaccessibility Scheme}

\author{ Shahram Mohsenipour}
\address{\ Shahram Mohsenipour,
         School of Mathematics, Institute for Research in Fundamental Sciences (IPM)\\
         P. O. Box 19395-5746, Tehran, Iran}
\email{mohseni@ipm.ir}
\thanks{This work was done while the author was a Postdoctoral Research Associate at the School of Mathematics,
        Institute for Research in Fundamental Sciences (IPM)}

\begin{abstract}
Suppose $L=\{<,\ldots\}$ is any countable first order language in which $<$ is interpreted as a linear order. Let $T$ be any complete first order theory in the language $L$ such that $T$ has a $\kappa$-like model where $\kappa$ is an inaccessible cardinal. Such $T$ satisfies the {\em Inaccessibility Scheme}. In this paper we study model theory of the inaccessibility scheme at the level of the existence of elementary end extensions for various models of it.
\end{abstract}

\maketitle
\bibliographystyle{amsplain}

\section{Introduction}

 First order schemes of sentences arising from set theoretical interests, can lead us to noticeable model theoretical investigations. The present paper can be considered as an attempt to proceed in this line. Let $L=\{<,\ldots\}$ be a countable first order language in which $<$ is interpreted as a linear order. Suppose $T$ is any complete first order theory in the language $L$ such that $T$ has a $\kappa$-like model $A$, where $\kappa$ is an inaccessible cardinal. The $\kappa$-likeness of $A$ means that $|A|=\kappa$ but for any $a\in A$, the cardinality of $\{x<a;x\in A\}$ is less than $\kappa$.
 It is easy to see that $T$ satisfies the {\em inaccessibility scheme:} It means that for any $\varphi(\bar{x},\bar{t})\in{L}$ we have

 \begin{center}
$(\forall z_{0}) (\exists z_{1}>z_{0})(\forall \bar{y}_{0})(\exists
\bar{y}_{1}<z_{1})E_{\varphi}(\bar{y}_{0}, \bar{y}_{1}; z_{0})$
\end{center}
\noindent where
\begin{center}
$E_{\varphi}(\bar{y}_{0}, \bar{y}_{1}; z_{0})= (\forall
\bar{t}<z_{0})(\varphi(\bar{y}_{0}, \bar{t})\equiv\varphi(\bar{y}_{1},
\bar{t})),$
\end{center}
note that in the above definition $\bar{t}$, $\bar{y}_{0}$ and $\bar{y}_{1}$ are tuples of variables with lengths clear from the context. In fact Shelah introduced the inaccessibility scheme in \cite{ShelahInaccScheme} for formulas $\varphi(\bar{x},\bar{t})\in{L}$ when $|\bar{x}|=1$, where $|\bar{x}|$ denotes the length of the tuple $\bar{x}$. But we will need to work with this more general form and call it again the inaccessibility scheme.

Now for the rest of the paper we fix a countable language $L=\{<,\ldots\}$ in which $<$ is interpreted as a linear order, a complete theory $T$ in the language $L$ with the Skolem functions.

Suppose $T$ satisfies the inaccessibility scheme. In this paper we are interested in the following question. Let $M$ be a model of $T$ with a given cardinality and cofinality. Does $M$ have an elementary end extension with a prescribed cardinality and cofinality? In general this is not true that any model of $T$ can have elementary end extensions of any possible cardinality. However it is interesting to know that if $T$ is any completion of PA, the Peano Arithmetic, then $T$ satisfies the above property. Though the above question is interesting in its own sake, our interest to it goes back mostly to the following problem, raised independently by Enayat \cite{aliunpublished} and Shelah \cite{shelahunderstand}: {\em Suppose $\lambda$ is an inaccessible but not a Mahlo cardinal. Does $T$ have a $\lambda$-like model?} See \cite{shahramsingular} for some historical backgrounds and some attempts on this problem and also see \cite{schmerltransfer} for some general backgrounds on the so called {\em model theoretic transfer theorems} and its relations with other branches of logic.

Before stating the results we recall some elementary concepts. For a regular cardinal $\kappa$ we say that a model $M$ of  $T$ has cofinality $\kappa$ and denote it by $cof(M)=\kappa$, if $\kappa$ is the least cardinal with the property that there is a strictly increasing cofinal sequence of elements of $M$ of length $\kappa$. Also for two models $M$ and $N$ of $T$ we say that $N$ is an elementary end extension of $M$ and denote it by $M\prec_{eee} N$, if $N$ is an elementary extension of $M$ and the order of $N$ end extends the order of $M$. Thorough the paper, $|M|$ denotes the cardinality of $M$. Shelah in \cite{ShelahInaccScheme} proved that

\begin{theorem}\label{0} Let $T$ satisfies the inaccessibility scheme. Suppose $M$ is a model of $T$ such that $M$ has cofinality $\omega$. Then $M$ has an elementary end extension $N$ such that $|M|=|N|$.
\end{theorem}

As our first theorem of this paper we will improve the above result by showing that $N$ can have cofinality $\omega$. More precisely:

\begin{theorem}\label{1} Let $T$ satisfies the inaccessibility scheme. Suppose $M$ is a model of $T$ such that $M$ has cofinality $\omega$. Then $M$ has an elementary end extension $N$ such that $|M|=|N|$ and $cof(N)=\omega$.
\end{theorem}

For our second theorem of this paper we will examine the above theorem for models of $T$ with uncountable cofinality and obtain the following theorem:

\begin{theorem}\label{2} Let $T$ satisfies the inaccessibility scheme. Suppose $M$ is a model of $T$ such that $M$ has cofinality $\kappa$ and is $\kappa$-saturated. Then for any regular cardinal $\lambda\leq\kappa$, $M$ has an elementary end extension $N$ such that $|M|=|N|$ and $cof(N)=\lambda$.
\end{theorem}

Our proofs of the both theorems rely heavily on Lemmas and Proposition \ref{3} about unboundedness of multivariables formulas which we will present and prove in the next section.

\section{Inaccessibility scheme and unbounded formulas}

\begin{definition} Suppose $\varphi(\bar{x},\bar{t})\in{L}$, where $\bar{x}=(x_1,\ldots,x_n)$. Let $M$ be a model of $T$ and $\bar{b}\in M$ such that $|\bar{b}|=|\bar{t}|$. We say that $\varphi(\bar{x},\bar{b})$ is \textsl{$\bar{x}$-unbounded} in $M$ if:
\begin{center}
$M\models(\forall\alpha_1)(\exists x_1>\alpha_1)\cdots(\forall\alpha_n)(\exists x_n>\alpha_n)\varphi(\bar{x},\bar{b}),$
\end{center}
otherwise we say that it is \textsl{$\bar{x}$-bounded}. We sometimes use the term ``unboundedness relative to $\bar{x}$" or just use ``unboundedness" when the set of variables is clear from the context.
\end{definition}
\begin{proposition}\label{3}
Let $T$ satisfies the inaccessibility scheme. Suppose $M$ is a model of $T$, $\psi(\bar{x},\bar{t})\in{L}$, $\bar{b}\in M$ and $|\bar{b}|=|\bar{t}|$ such that $\psi(\bar{x},\bar{b})$ is $\bar{x}$-unbounded in $M$. Suppose $\varphi(\bar{x},\bar{t})\in{L}$ and $a\in M$ are given. Then there is a tuple $\bar{c}\in M$ with $|\bar{c}|=|\bar{x}|$ such that $E_{\varphi}(\bar{x},\bar{c};a)$$\wedge$$\psi(\bar{x},\bar{b})$ is $\bar{x}$-unbounded in $M$.
\end{proposition}
\begin{proof}
The proof proceeds by induction on the length of $|\bar{x}|=n$. The case $n=1$ is due to Shelah \cite{ShelahInaccScheme}. Suppose $n>1$. By the inaccessibility scheme there is $y_1$ such that for all $\bar{x}$ there exists $\bar{d}<y_1$ such that
\begin{center}
$M\models E_{\varphi}(\bar{x},\bar{d};a)$.
\end{center}
Suppose for each $\bar{d}<y_1$ the formula $E_{\varphi}(\bar{x},\bar{d};a)\wedge\psi(\bar{x},\bar{b})$ is $\bar{x}$-bounded. So by the definition of unboundedness we have:
\begin{center}
$M\models(\exists\alpha_1)(\forall x_1>\alpha_1)\cdots(\exists\alpha_n)(\forall x_n>\alpha_n)\neg[E_{\varphi}(\bar{x},\bar{d};a)\wedge\psi(\bar{x},\bar{b})].$
\end{center}
Now we want to detect $\alpha_1$ by using a suitable Skolem function. So suppose $\theta(h,\bar{u},v,\bar{w})$ is the following formula:
\begin{center}
$(\forall x_1>h)(\exists\alpha_2)(\forall x_2>\alpha_2)\cdots(\exists\alpha_n)(\forall x_n>\alpha_n)\neg[E_{\varphi}(\bar{x},\bar{u};v)\wedge\psi(\bar{x},\bar{w})]$,
\end{center}
where the length of the tuples are clear from the context. Let $\tau(\bar{u},v,\bar{w})=h$ is the Skolem function of $\theta(h,\bar{u},v,\bar{w})$. Obviously $M\models(\exists\alpha)\theta(\alpha,\bar{d},a,\bar{b})$. Also $\tau(\bar{u},a,\bar{b})=\alpha$ implies $M\models(\exists\alpha)\theta(\alpha,\bar{u},a,\bar{b})$. We claim that the set $\Gamma=\{\alpha\in M: (\exists\bar{u}<y_1)\tau(\bar{u},a,\bar{b})=\alpha\}$ is unbounded in $M$. Suppose it is not the case and $A\in M$ is an upper bound of $\Gamma$. Fix an $x^{*}>A$, so for each $\bar{d}<y_1$ we have:

\vspace{.1in}\hspace{.2in}$M\models(\exists\alpha_2)(\forall x_2>\alpha_2)\cdots(\exists\alpha_n)(\forall x_n>\alpha_n)$\vspace{.07in}

\hspace{1.5in}$\neg[E_{\varphi}(x^{*},x_2,\ldots,x_n,\bar{d};a)\wedge\psi(x^{*},x_2,\ldots,x_n,\bar{b})].$\vspace{.1in}

\noindent Then we conclude that $E_{\varphi}(x^{*},x_2,\ldots,x_n,\bar{d};a)\wedge\psi(x^{*},x_2,\ldots,x_n,\bar{b})$ is $(x_2,\ldots,x_n)$-bounded in $M$ ($\clubsuit$). Clearly $\psi(x^{*},x_2,\ldots,x_n,\bar{b})$ is  $(x_2,\ldots,x_n)$-unbounded in $M$, then by the induction hypothesis there is a tuple $\bar{g}$ in $M$ with $|\bar{g}|=n-1$ such that $E_{\varphi}(x^{*},x_2,\ldots,x_n,x^{*},\bar{g};a)\wedge\psi(x^{*},x_2,\ldots,x_n,\bar{b})$ is $(x_2,\ldots,x_n)$-unbounded in $M$. Again by using the inaccessibility scheme there is a tuple $\bar{e}<y_1$ (for $y_1$ recall the first paragraph of the proof), such that
\begin{center}
$M\models\forall\bar{t}<a[\varphi(\bar{e},\bar{t})\equiv\varphi(x^{*},\bar{g},\bar{t})]$
\end{center}
which immediately implies that $E_{\varphi}(x^{*},x_2,\ldots,x_n,\bar{e};a)$$\wedge$$\psi(x^{*},x_2,\ldots,x_n,\bar{b})$ is unbounded relative to $(x_2,\ldots,x_n)$, which is in contradiction with ($\clubsuit$). So we have proved the claim that $\Gamma$ is unbounded. To finish the proof recall our Skolem function $\tau(\bar{u},v,\bar{w})=h$. Suppose $f>a,\bar{b},y_1$, by another use of the inaccessibility scheme there exists $y_2$ such that for all $h$ there exists $h^{'}<y_2$ such that
\begin{center}
$M\models\forall\bar{t},v,\bar{w}<f[\tau(\bar{t},v,\bar{w})=h \equiv \tau(\bar{t},v,\bar{w})=h^{'}]$
\end{center}
which implies that for any $h>y_2$ we have $\forall\bar{t}<y_1\tau(\bar{t},a,\bar{b})\neq h$. So $y_2$ is an upper bound of $\Gamma$, which violates the claim. Therefore the proposition is proved.
\end{proof}

We also will need another scheme of first order sentences in the course of our proofs, which is called the {\em Regularity Scheme}. It is introduced by Keisler in \cite{keislertree} and says: {\em Suppose $\varphi(\bar{x},\bar{y},\bar{t})$ is any formula in $L$. If $\{\bar{x};(\exists\bar{y}<y_0)\varphi(\bar{x},\bar{y},\bar{t})\}$ is $\bar{x}$-unbounded, then for some $\bar{y_1}<y_0$, the set $\{\bar{x};\varphi(\bar{x},\bar{y_1},\bar{t})\}$ is $\bar{x}$-unbounded.} In fact Keisler considered the regularity scheme for one-variable formulas, namely, the case with $|\bar{x}|=1$. But as in the inaccessibility case, we prefer the general form. It is easy to see that if $T$ has a model with an elementary end extension then $T$ satisfies the regularity scheme for the case $|\bar{x}|=1$. By the Shelah Theorem \ref{0} if $T$ satisfies the inaccessibility scheme then every countable model of $T$ has an elementary end extension, therefore the inaccessibility scheme implies the regularity scheme for the case $|\bar{x}|=1$. But it seems this argument does not work for the general case of the regularity scheme. In the proof of the following lemma we show a direct proof of the above implication and also show that the regularity scheme for formulas with $|\bar{x}|=1$ implies the general regularity scheme:

\begin{lemma}\label{4}
(i)Suppose $T$ satisfies the inaccessibility scheme, then $T$ satisfies the regularity scheme for formulas with $|\bar{x}|=1$.

(ii)Suppose $T$ satisfies the regularity scheme for formulas with $|\bar{x}|=1$, then $T$ satisfies the regularity scheme.
\end{lemma}
\begin{proof}
(i) Let $\varphi(x,\bar{y},\bar{t})$ be a formula and $M$ be a model of $T$. Suppose $y_0,\bar{b}\in M$ with $|\bar{t}|=|\bar{b}|$ such that for every $\bar{y}<y_0$ the set $A=\{x\in M;\varphi(x,\bar{y},\bar{b})\}$ is $x$-bounded in $M$. It suffices to show that $B=\{x\in M;(\exists\bar{y}<y_0)\varphi(x,\bar{y},\bar{b})\}$ is $x$-bounded in $M$. $A$ is being bounded implies that
\begin{center}
$M\models (\forall\bar{y}<y_0)(\exists\alpha)(\forall x>\alpha)\neg\varphi(x,\bar{y},\bar{b})$.
\end{center}
Now consider the related Skolem function $\tau(\bar{y},\bar{b})=\alpha$ and set
\begin{center}
$\Gamma=\{\alpha\in M;(\exists\bar{y}<y_{0})\tau(\bar{y},\bar{b})=\alpha\}.$
\end{center}
We show that $\Gamma$ is bounded in $M$. By the inaccessibility scheme
\begin{center}
$M\models(\exists y_1)(\forall\alpha)(\exists\alpha^{'}<y_1)(\forall\bar{y}<y_0)[\tau(\bar{y},\bar{b})=\alpha \equiv \tau(\bar{y},\bar{b})=\alpha^{'}].$
\end{center}
Therefore if $\alpha>y_1$, then $(\forall\bar{y}<y_0)\tau(\bar{y},\bar{b})\neq\alpha$ which establishes that $y_1$ is an upper bound of $\Gamma$. Recalling the definition of the Skolem function $\tau(\bar{y},\bar{b})=\alpha$, we deduce that
\begin{center}
$M\models(\forall x >y_1)(\forall\bar{y}<y_0)\neg\varphi(x,\bar{y},\bar{b})$.
\end{center}
Thus $y_1$ is an upper bound of $B$. This proves (i).

(ii) The proof is by induction on $|\bar{x}|=n$. Suppose $n>1$. Let $\varphi(\bar{x},\bar{y},\bar{t})$ be a formula with $\bar{x}=(x_1,\ldots,x_n)$ and $M$ be a model of $T$. Suppose $y_0,\bar{b}\in M$ with $|\bar{t}|=|\bar{b}|$ such that for every $\bar{y}<y_0$ the set $A=\{\bar{x}\in M;\varphi(\bar{x},\bar{y},\bar{b})\}$ is $\bar{x}$-bounded in $M$. More formally:
\begin{center}
$M\models(\forall\bar{y}<y_0)(\exists\alpha_1)(\forall x_1>\alpha_1)\cdots(\exists\alpha_n)(\forall x_n>\alpha_n)\neg\varphi(\bar{x},\bar{y},\bar{b})$.
\end{center}
Now consider the Skolem function $\tau(\bar{y},\bar{b})=\alpha_1$ of the formula
\begin{center}
$(\forall x_1>\alpha_1)(\exists\alpha_2)(\forall x_2>\alpha_2)\cdots(\exists\alpha_n)(\forall x_n>\alpha_n)\neg\varphi(\bar{x},\bar{y},\bar{b})$,
\end{center}
by the regularity scheme for one-variable formulas we deduce that the set
\begin{center}
$\{\alpha\in M;(\exists\bar{y}<y_0)\tau(\bar{y},\bar{b})=\alpha\}$
\end{center}
is bounded in $M$. Let $A$ be an upper bound and fix an $x^{*}>A$. This means that for all $\bar{y}<y_0$ the set $\{(x_2,\ldots,x_n)\in M;M\models\varphi(x^{*},x_2,\ldots,x_n,\bar{y},\bar{b})\}$ is $(x_2,\ldots,x_n)$-bounded in $M$. Now by the induction hypothesis we conclude that the set
\begin{center}
$\{(x_2,\ldots,x_n)\in M;M\models(\exists\bar{y}<y_0)\varphi(x^{*},x_2,\ldots,x_n,\bar{y},\bar{b})\}$
\end{center}
is $(x_2,\ldots,x_n)$-bounded in $M$. Since $x^{*}$ can be arbitrary large, it follows that the set $\{\bar{x};M\models(\exists\bar{y}<y_0)\varphi(\bar{x},\bar{y},\bar{b})\}$ is $\bar{x}$-bounded. This completes the proof of (ii).
\end{proof}

The next lemma will enable us to control the cofinality of elementary end extensions. Suppose $x_1,\ldots,x_n$ are any variables and $\tau(\bar{u},\bar{t})$ is a term such that $k=|\bar{u}|<n$. For any integers $0<i_1<\cdots<i_{k+1}\leq n$ and any $a$ we put
\begin{center}
$A_{\tau}(x_{i_1},\ldots,x_{i_{k+1}};a)=\forall\bar{t}<a[\tau(x_{i_1},\ldots,x_{i_{k}},\bar{t})<x_{i_{k+1}}]$.
\end{center}
and
\begin{center}
$B_{\tau}(x_1,\ldots,x_n;a)=\bigwedge_{0<i_1<\cdots<i_{k+1}\leq n}A_{\tau}(x_{i_1},\ldots,x_{i_{k+1}};a)$
\end{center}
Clearly the above definition of $B_{\tau}$ depends on $n$, but for simplicity we don't exhibit $n$. Now we can state the next lemma:
\begin{lemma}\label{5}
Let $T$ satisfies the regularity scheme. Suppose $M$ is a model of $T$, $\varphi(\bar{x},\bar{t})\in{L}$, $\bar{b}\in M$ and $|\bar{b}|=|\bar{t}|$ such that $\varphi(\bar{x},\bar{b})$ is $\bar{x}$-unbounded in $M$. Suppose $\tau(\bar{u},\bar{t})\in{L}$ is a term such that $|\bar{u}|<|\bar{x}|$ and $a\in M$. Then $\varphi(\bar{x},\bar{b})\wedge B_{\tau}(\bar{x};a)$ is $\bar{x}$-unbounded in $M$.
\end{lemma}
\begin{proof}
It suffices to show that for any $1\leq i_1<\cdots<i_{k+1}\leq n$,
\begin{center}
$A_{\tau}(x_{i_1},\ldots,x_{i_{k+1}};a)\wedge\varphi(\bar{x},\bar{b})$
\end{center}
is  $\bar{x}$-unbounded in $M$. The rest of the proof is better explainable in a simple game theoretical language. For any formula $\psi(x_1,\ldots,x_n,\bar{c})$ where $\bar{c}$ is from $M$, we introduce a game of length $n$, $G(\psi(x_1,\ldots,x_n,\bar{c}),M)$, for two players I and II such that each player plays $n$ moves. Player I chooses for his $i$th move an element $\alpha_i$ in $M$ and then the player II chooses for her $i$th move an element $x_i>\alpha_i$ in $M$. Player I wins the game if $M\models\neg\psi(x_1,\ldots,x_n,\bar{c})$, otherwise the player II wins the game. It is easy to see that $\psi(x_1,\ldots,x_n,\bar{c})$ is $\bar{x}$-unbounded in $M$ iff the player II has a winning strategy. By the assumption we know that the player II has a winning strategy for the game $G_1=G(\varphi(\bar{x},\bar{b}),M)$. Now we show that she has a winning strategy too for the game
$G_2=G(A_{\tau}(x_{i_1},\ldots,x_{i_{k+1}};a)\wedge\varphi(\bar{x},\bar{b}),M)$. For any $i<i_{k+1}$ the player II plays the $i$th move according to her winning strategy in the game $G_1$. For her $i_{k+1}$th move she checks the previously played moves $x_{i_1},\ldots,x_{i_{k}}$, but we know that by the regularity scheme, the set $\{y\in M;(\exists\bar{t}<a)\tau(x_{i_1},\ldots,x_{i_{k}},\bar{t})=y\}$ is bounded in $M$, say by $\beta$, then the player II chooses for her $i_{k+1}$th move any element greater than max$\{x_{i_{k+1}},\beta\}$ in which $x_{i_{k+1}}$ is her $i_{k+1}$th move for the game $G_1$.  For the rest of the game, the player II plays according to her winning strategy for the game $G_1$. It is not hard to see that this is really a winning strategy for the player II in the game $G_2$.
\end{proof}
\begin{remark}\label{rem}
In Proposition \ref{3} and Lemma \ref{5}, we assumed that the two formulas that their conjunction were unbounded, had the same set of variables. It is not necessary. Just notice that by adding trivial formulas $(x_i=x_i)$, if necessary, we can lie in the situation of the above Proposition and Lemma.

\end{remark}

\section{Proofs of the main theorems}
Now we are ready to prove our stated theorems in the Introduction.
\begin{proof}
[Proof of Theorem \ref{1}] Let $\{a_i;i\in\omega\}$ is a cofinal set of elements of $M$. Fix a countable set of variables $X=\{x_1,x_2,\ldots,x_n,\ldots\}$. Suppose $\{\varphi_i(\bar{x}_i,\bar{t}_i);i\in\omega\}$ and $\{\tau_i(\bar{y}_i,\bar{s}_i);i\in\omega\}$ enumerates all formulas and terms of $L$ such that $\bar{x}_i$'s and $\bar{y}_i$'s cover all finite sequences of elements of $X$ (of course in the increasing order of indices) and each formula and each term occurs infinitely often. We inductively construct a type $p(x_1,x_2,\ldots)$ consisting of formulas with parameters in $M$ which is finitely satisfiable in $M$. By Proposition \ref{3}, Lemma \ref{5} and Remark \ref{rem}, there is a tuple $\bar{c}_0\in M$ with the same length as $\bar{x}_0$, such that
\begin{center}
$E_{\varphi_0}(\bar{x}_0,\bar{c}_0;a_0)\wedge B_{\tau_0}(\bar{y}_0;a_0)$
\end{center}
is $(\bar{x}_0,\bar{y}_0)$-unbounded in $M$. We denote the above formula by $A_0(\bar{x}_0,\bar{y}_0)$. Suppose $A_n(\bar{x}_0,\bar{y}_0,\ldots,\bar{x}_n,\bar{y}_n)$ has been defined and is $(\bar{x}_0,\bar{y}_0,\ldots,\bar{x}_n,\bar{y}_n)$-unbounded in $M$. Again by Proposition \ref{3}, Lemma \ref{5} and Remark \ref{rem}, there is a tuple $\bar{c}_{n+1}\in M$ with the same length as $\bar{x}_{n+1}$, such that
\begin{center}
$A_n(\bar{x}_0,\bar{y}_0,\ldots,\bar{x}_n,\bar{y}_n)\wedge E_{\varphi_{n+1}}(\bar{x}_{n+1},\bar{c}_{n+1};a_{n+1})\wedge B_{\tau_{n+1}}(\bar{y}_{n+1};a_{n+1})$
\end{center}
is $(\bar{x}_0,\bar{y}_0,\ldots,\bar{x}_{n+1},\bar{y}_{n+1})$-unbounded in $M$. We denote the above formula by $A_{n+1}(\bar{x}_0,\bar{y}_0,\ldots,\bar{x}_{n+1},\bar{y}_{n+1})$ and put
\begin{center}
$p(x_1,x_2,\ldots)=\{A_i(\bar{x}_0,\bar{y}_0,\ldots,\bar{x}_i,\bar{y}_i);i\in\omega\}$.
\end{center}
Clearly the above type is finitely satisfiable in $M$, so there is an elementary extension $N^{*}$ of $M$ such that $N^{*}$ realizes $p(x_1,x_2,\ldots)$. Suppose $X^{*}=\{x_1^{*},x_2^{*},\ldots,x_n^{*},\ldots\}$ is a set of realizations in $N^{*}$. Form the Skolem hull of $M\cup X^{*}$ in $N^{*}$ and denote it by $N$. We will show that $N$ is an elementary end extension of $M$ with $|N|=|M|$ and $cof(N)=\omega$. Since the language $L$ is countable it is clear that $|N|=|M|$. By construction, for every term $\tau(x_1,x_2,\ldots,x_n,\bar{t})$ and $a\in M$ and every integers $1\leq i_1<\cdots<i_{n+1}$, we have
\begin{center}
$N\models\forall\bar{t}<a[\tau(x_{i_1}^{*},\ldots,x_{i_{n}}^{*},\bar{t})<x_{i_{n+1}}^{*}]$.
\end{center}
Thus $X^{*}=\{x_1^{*},x_2^{*},\ldots,x_n^{*},\ldots\}$ is cofinal in $N$. To show $M\prec_{eee} N$, suppose for some $x_{i_1}^{*},\ldots,x_{i_{n}}^{*}$ and $\bar{b}\in M$, $N\models\tau(x_{i_1}^{*},\ldots,x_{i_{n}}^{*},\bar{b})<a_k$ for some $k\in\omega$ and term $\tau(\bar{u},\bar{t})\in L$. By recalling the process of finding tuples $\bar{c}_i\in M$ in constructing the formulas $E_{\varphi_i}(\bar{x}_i,\bar{c}_i;a_i)$'s and also noting that each formula $\varphi_i(\bar{x}_i,\bar{t}_i)$ appears infinitely often in the course of the construction, we deduce that there are $m>k$ such that $a_{m+1}>\bar{b}$ and $c_{i_1},\ldots,c_{i_n}$ in $M$ so that:
\begin{equation}
N\models(\forall\bar{t}<a_{m+1})[\tau(c_{i_1},\ldots,c_{i_n},\bar{t})<a_m \equiv \tau(x_{i_1}^{*},\ldots,x_{i_{n}}^{*},\bar{t})<a_m)],
\end{equation}
as well as
\begin{equation}
N\models(\forall\bar{t},y<a_{m+1})[\tau(c_{i_1},\ldots,c_{i_n},\bar{t})=y \equiv \tau(x_{i_1}^{*},\ldots,x_{i_{n}}^{*},\bar{t})=y)].
\end{equation}
Suppose $\tau(c_{i_1},\ldots,c_{i_n},\bar{b})=d\in M$, by (1) we deduce that $N\models d<a_m$ and by (2) we conclude that $\tau(x_{i_1}^{*},\ldots,x_{i_{n}}^{*},\bar{b})=d\in M$. Thus $N$ is an elementary end extension of $M$. This finishes the proof of Theorem \ref{1}.
\end{proof}

By iterating Theorem \ref{1} $\omega_1$ times, we immediately obtain the following corollary:

\begin{corollary}
 Let $T$ satisfies the inaccessibility scheme. Suppose $M$ is a model of $T$ such that $M$ has cofinality $\omega$. Then $M$ has an elementary end extension $N$ such that $cof(N)=\omega_{1}$. If $|M|$ is uncountable then $|M|=|N|.$
\end{corollary}

Now we turn to prove Theorem \ref{2}. In most parts the proof is the same as the proof of Theorem \ref{1}.

\begin{proof}
[Proof of Theorem \ref{2}] We start similar to the proof of Theorem \ref{1}. Let $\{a_i;i<\kappa\}$ is a cofinal set of elements of $M$. Fix a set of variables $X=\{x_1,\ldots,x_i,\ldots;i<\lambda\}$. Suppose $\{\varphi_i(\bar{x}_i,\bar{t}_i);i<\lambda\}$ and $\{\tau_i(\bar{y}_i,\bar{s}_i);i<\lambda\}$ enumerates all formulas and terms of $L$ such that $\bar{x}_i$'s and $\bar{y}_i$'s cover all finite sequences of elements of $X$ (of course in the increasing order of indices) and each formula and each term occurs $\kappa$ times. We inductively construct a type $p(x_1,\ldots,x_i,\ldots)_{i<\lambda}$ consisting of formulas with parameters in $M$ which is finitely satisfiable in $M$. We  construct for each $i<\omega$, the formulas $A_i(\bar{x}_0,\bar{y}_0,\ldots,\bar{x}_i,\bar{y}_i)$ as in the proof of Theorem \ref{1}. For simplicity, for $i<\omega$, we put $\bar{X}_i=(\bar{x}_0,\bar{y}_0,\ldots,\bar{x}_i,\bar{y}_i)$. Now for any $\omega\leq i<\kappa$, we define a formula $A_i(\bar{X}_i)$ with the parameters from $M$, where $\bar{X}_i=(\bar{x}_i,\bar{y}_i)$ such that
\begin{center}
$p(x_1,\ldots,x_i,\ldots)_{i<\lambda}=\{A_i(\bar{X}_i);i<\kappa\}$.
\end{center}
We construct $A_{\omega}(\bar{X}_{\omega})$, the construction of all $A_i(\bar{X}_i)$ for $\omega< i<\kappa$ will be similar. For any finite subset $J\subset\omega$, by Proposition \ref{3}, Lemma \ref{5} and Remark \ref{rem}, we deduce that there exists a tuple $\bar{v}\in M$ with $|\bar{v}|=|\bar{x}_{\omega}|$ such that the formula
\begin{equation} \label{A}
\bigwedge_{i\in J}A_i(\bar{X}_i)\wedge E_{\varphi_{\omega}}(\bar{x}_{\omega},\bar{v},a_{\omega})
\wedge B_{\tau_{\omega}}(\bar{y}_{\omega};a_{\omega})
\end{equation}
is $(\bar{X}_i,\bar{x}_{\omega},\bar{y}_{\omega})_{i\in J}$-unbounded in $M$. Since $M$ is $\kappa$-saturated, then there exists a tuple  $\bar{c}_{\omega}\in M$ such that it can be served as $\bar{v}$ for the formula (\ref{A}) for all finite $J\subset{\omega}$. Now we put
\begin{center}
$A_{\omega}(\bar{X}_{\omega})=E_{\varphi_{\omega}}(\bar{x}_{\omega},\bar{c}_{\omega},a_{\omega})
\wedge B_{\tau_{\omega}}(\bar{y}_{\omega};a_{\omega})$
\end{center}
We can define in the same manner all $A_i(\bar{X}_i)$ for $\omega<i<\kappa$. For example to define $A_{\omega+1}(\bar{X}_{\omega+1})$, we must consider all finite subsets $J\subset \omega+1$. Now the definition of the type $p(x_1,\ldots,x_i,\ldots)_{i<\lambda}$ is complete and it is finitely satisfiable in $M$. Therefore there is an elementary extension $N^{*}$ of $M$ such that $N^{*}$ realizes $p(x_1,\ldots,x_i,\ldots)_{i<\lambda}$. Suppose $X^{*}=\{x_1^{*},\ldots,x_i^{*},\ldots;i<\lambda\}$ is a set of realizations in $N^{*}$. Form the Skolem hull of $M\cup X^{*}$ in $N^{*}$ and denote it by $N$. By reasoning exactly in the same way as the proof of Theorem \ref{1}, we can show that $X^{*}=\{x_1^{*},\ldots,x_i^{*},\ldots;i<\lambda\}$ is cofinal in $N$ and $M\prec_{eee}N$.
\end{proof}

\section{Concluding remarks and open questions}
In a previous paper \cite{shahramkeislermorley}, we introduced another first order scheme of axioms, called the Erd\"{o}s-Rado scheme. We now introduce a more general and more suitable version: For all formulas $\varphi(\bar{x},\bar{t})$ and all terms $\tau(\bar{y},\bar{s})$ such that $|\bar{y}|<|\bar{x}|$ and also all $m\geq|\bar{x}|$ we have: For all $z_0$ and all $z_1>z_0$, there exist $y_1<\cdots<y_m$ greater than $z_1$ such that:
\begin{center}
$\bigwedge_{\bar{u},\bar{v}\subset\bar{y}}E_{\varphi}(\bar{u},\bar{v};z_0)\wedge B_{\tau}(y_1,\ldots,y_m;z_0).$
\end{center}
By using the Erd\"{o}s-Rado partition theorem, it is easily seen that if $T$ has a $\kappa$-like model where $\kappa$ is an inaccessible cardinal then $T$ satisfies the Erd\"{o}s-Rado scheme. In \cite{shahramkeislermorley} we proved a theorem which its proof works as well for this new version of the Erd\"{o}s-Rado scheme:

\begin{theorem}\label{6}  Suppose $T$ satisfies the Erd\"{o}s-Rado scheme and the inaccessibility scheme and $M$ is a countable model of $T$ which is $\omega$-saturated. Then for any regular cardinal $\kappa$ and any cardinal $\lambda\geq\kappa$, $M$ has an elementary end extension with cardinality $\lambda$ and cofinality $\kappa$.
\end{theorem}

By carefully checking the proof and making easy changes, we can eliminate the countability assumption and replace the inaccessibility scheme by the regularity scheme. More precisely:

\begin{theorem}\label{7}  Suppose $T$ satisfies the Erd\"{o}s-Rado scheme and the regularity scheme and $M$ is model of $T$ which is $\omega$-saturated and $cof(M)=\omega$. Then for any regular cardinal $\kappa$ and any cardinal $\lambda\geq\kappa$ such that $\lambda\geq|M|$, $M$ has an elementary end extension with cardinality $\lambda$ and cofinality $\kappa$.
\end{theorem}

Two questions come to mind:
\begin{question}
What is the exact relation between the Erd\"{o}s-Rado scheme and the inaccessibility scheme? Does one of them imply the other?
\end{question}
We would like to add to the above question the {\em reflection scheme} too, which is studied in \cite{shahramali}.

Regarding Enayat and Shelah's open question, mentioned in the introduction, we ask:
\begin{question}\label{Q}
Suppose $\lambda$ is a singular cardinal and $T$ satisfies the Erd\"{o}s-Rado scheme and the inaccessibility scheme. Does $T$ have a $\lambda$-like model?
\end{question}
Keisler has proved in \cite{keislerordering} that if $T$ has a $\kappa$-like model when $\kappa$ is a strong limit cardinal, then $T$ has a $\lambda$-like model, when $\lambda$ is a singular cardinal. So positive answer to Question \ref{Q} will contrast with Keisler's theorem. To answer the above question, one might try to formalize Keisler's proof in $T$.

\bibliography{inaccessibility}
\bibliographystyle{plain}
\end{document}